\documentclass[10pt]{amsart}

\usepackage{amsmath}
\usepackage{amsthm}
\usepackage{amsfonts}
\usepackage{amssymb,latexsym}
\usepackage[all]{xy}
\usepackage{graphicx}
\usepackage[mathscr]{eucal}

\theoremstyle{plain}
    \newtheorem{thm}{Theorem}

\newcommand{\C}{\mathcal C}
\newcommand{\D}{\mathcal D}

\newcommand{\F}{\mathcal F}

\newcommand{\la}{\ensuremath{\langle}}
\newcommand{\ra}{\ensuremath{\rangle}}
\newcommand{\rar}{\ensuremath{\rightarrow}}

\begin{document}

%
\title[Thick subcategories in stable homotopy theory]{Thick subcategories in stable homotopy theory \\
(work of Devinatz, Hopkins, and Smith).}
\date{\today}

\author{Sunil K. Chebolu}

\maketitle

In this series of lectures we give an exposition of the seminal work
of Devinatz, Hopkins, and Smith  which is surrounding the classification of the thick subcategories of finite
spectra in stable homotopy theory. The lectures are expository and are aimed primarily at
non-homotopy theorists. We begin with an introduction to the stable
homotopy category of spectra, and then talk about the celebrated
thick subcategory theorem and discuss a few  applications to the
structure of the Bousfield lattice. Most of the results that we
discuss below were conjectured by Ravenel \cite{rav} and were proved
by Devinatz, Hopkins, and Smith \cite{dhs, hs}.

\section{The stable homotopy category of spectra}
Recall that in homotopy theory one is interested in studying the
homotopy classes of maps between CW complexes  (spaces that are
built in a systematic way by attaching cells): If $f$ and $g$ are
maps (continuous)  between CW complexes $X$ and $Y$, we say that
they are  \emph{homotopic} if there is a map from the cylinder $X
\times [0, 1]$ to $Y$ whose restriction to the two  ends (top and
bottom) of the cylinder gives $f$ and $g$ respectively. The homotopy
classes of maps between $X$ and $Y$ is denoted by $[X, Y]$.
In stable homotopy theory one studies a
weaker notion of homotopy called stable homotopy
 -- maps $f$ and $g$ as above are said to be \emph{stably homotopic} if $\Sigma^n f$ and
$\Sigma^n g$ are homotopic for some $n$. ($\Sigma$ denotes the
reduced suspension functor on the homotopy category of pointed CW
complexes.)  The notion of stable homotopy is much weaker than
homotopy. For example, the obvious quotient map from the torus to
the two sphere is not null homotopic but stably null homotopic.
The importance of  stable homotopy classes of maps comes from an old
result due to Freudenthal which implies that if $X$ and $Y$ are
finite CW complexes, then the sequence
\[ [X, Y] \rar [\Sigma X, \Sigma Y] \rar [\Sigma^2 X, \Sigma^2 Y] \rar \cdots \]
eventually stabilises. The stable homotopy classes of maps from $X$
to $Y$ is precisely the above colimit. In particular when $X$ is the
$n$-sphere $S^n$, we get the $n$-th stable homotopy group of $Y$,
denoted $\pi^s_n(Y)$. Computing stable homotopy groups is, in general, a more manageable problem than
that of homotopy groups. However, it became abundantly clear to homotopy theorists by 1960s that in order to do serious
stable calculations efficiently it is absolutely essential to have a nice
category in which the objects are stabilised analogue of spaces each of which represent  a cohomology theory.  The finite objects of such a
category can be easily described. This is called the (finite)
Spanier-Whitehead category which captures finite stable phenomena,
and is defined as follows. The objects are ordered pairs $(X, n)$
where $X$ is a finite CW complex and $n$ is an integer, and
morphisms between objects $(X ,n)$ and $(Y, m)$ are given by
\[\{(X,n), (Y,m)\}:= \underset{k}{\text{colim}}\; [\Sigma^{n+k} X , \Sigma^{m+k} Y ].\]
This category has a formal suspension $\Sigma (X, n) := (X, n+1)$
which agrees with the geometric suspension, i.e., $(X , n+1) \cong
(\Sigma X , n)$. While there is no geometric desuspension, there is
a formal desuspension $\Sigma^{-1} (X, n) = (X , n-1)$. Thus by passing to the Spanier-Whitehead category we have
inverted the suspension functor on CW complexes! Moreover, this
category  has a tensor triangulated structure: exact
triangles are induced by mapping sequences and the product comes from the
smash product of CW complexes. Although this category is the right stabilisation of finite CW complexes, it has its limitations. The key point here is that
one needs infinite dimensional CW complexes to understand finite CW complexes.
For example, the singular cohomology theories on finite CW complexes are represented by
Eilenberg-Mac Lane spaces which are infinite dimensional. So  naturally one has to
enlarge the finite Spanier-Whitehead category so that it has all the
desired properties; one at least demands arbitrary coproducts and
the Brown representability theorem. Building the ``stable category''
with all the desired properties is quite challenging. Several stable categories have been proposed; the first satisfactory category was constructed by
Mike Boardman in his 1964 Warwick thesis \cite{boardman-thesis}, and then by Frank Adams \cite{bluebook}, followed by several others.
All these models share a set of properties which can be taken to be the defining properties of the stable homotopy
category. Following Margolis \cite{mar}, we take this axiomatic
approach.

\begin{thm} There is a category $\mathcal{S}$ called the stable
homotopy category (whose objects are called spectra)
which has the following properties.

\begin{enumerate}
\item $\mathcal{S}$ is a triangulated category which admits arbitrary set indexed coproducts.
\item $\mathcal{S}$ has a unital, commutative and associative smash product which is compatible with
the triangulation.
\item The sphere spectrum is a graded weak generator: $\pi_*(X) = 0$
implies $X = 0$.
\item The full subcategory of compact objects of $\mathcal{S}$ is
equivalent to the Spanier-Whitehead category of finite CW complexes.
\end{enumerate}
\end{thm}

Note that there are a lot of categories in the literature which satisfy the first three properties. It is the
property (4) that makes the theorem very unique, important and non-trivial.
It is also worth pointing out that the study of spectra is
equivalent to that of generalised homology theories on CW complexes
(theories which satisfy all the Eilenberg-Steenrod axioms except the
dimension axiom.) Some standard examples of such theories are the
singular homology, complex K-theory are Complex bordism which are
represented by the Eilenberg-Mac Lane spectrum, the K-theory spectrum
$BU$ and the Thom spectrum $MU$ respectively. The study of these two subjects is in turn essentially
equivalent to the study of infinite loop spaces.
To get a better picture of the strong connections between spectra,
generalised homology theories and infinite loop spaces,  we refer the reader to Adams excellent account \cite{adams-infinite}.

The stable homotopy category is very rich in its structural complexity, and one of
the goals of the subject is to understand the
global structure of this category. Doug Ravenel in the late 70s
suspected some deep and interesting structure in this category
(which was inspired by his algebraic calculations) and has
formulated seven conjectures \cite{rav} on the structure
of $\mathcal{S}$. All but one of them have been solved by 1986, due largely
to the seminal work of Devinatz, Hopkins, and Smith \cite{dhs, hs}. We discuss some of these conjectures which are
surrounding the thick subcategory theorem.

To start, let $f: X \rar Y$ be a map between spectra. Then we can ask
several questions, the first one is when is $f$ null-homotopic?
Detecting null homotopy of maps is an extremely difficult problem. A
long standing conjecture of Peter Freyd \cite{freydGH} called Generating Hypothesis says that if $X$ and $Y$ are
finite spectra, then $f$ is null homotopic if $\pi_*(f)$ is zero.
Some partial results are known due to Devinatz \cite{ethanGH}, the conjecture remains
open; see \cite{freydGH} for some very interesting consequences of this
conjecture. The second question is when is $f$ nilpotent under
composition. The nilpotence theorem which was conjectured by Ravenel
gives an answer to this question when the spectra in question are
finite. This theorem is very deep and its proof involves some hard homotopy
theory. It generalises a well-known theorem of Nishida which tells that every
positive degree self map of the sphere spectrum is nilpotent.

\begin{thm}\cite{dhs} \emph{(Nilpotence theorem)} There is a generalised homology theory known as
$MU_*(-)$ (complex bordism) such that a map $f: X \rar Y$ between finite spectra is
nilpotent if and only if $MU_*(f)$ is nilpotent.
\end{thm}

A much sharper view of the stable homotopy theory is obtained when
one localises at the prime $p$ and studies the $p$-local stable
category whose objects are spectra whose homotopy groups are
$p$-local, i.e., $\pi_*(X) \cong \pi_*(X) \otimes \mathbb{Z}_{(p)}.$
It is a very standard practise in stable homotopy theory to localise
at a prime $p$. When this is done, there are distinguished field
objects known as Morava K-theories $K(n)$ (with the prime $p$
suppressed) which play a key role in  the $p$-local stable
category. We now begin discussing these objects which also play an important
role in the thick subcategory theorem.

\section{Morava K-theories and the thick subcategory theorem.}

To set the stage, let $\F$ denote the category of compact objects in
the p-local stable homotopy category $\mathcal{S}$. There are many
naturally arising properties of spectra called \emph{generic properties} which are properties that are preserved
under cofibrations, retractions and suspensions. Recall that a
subcategory is \emph{thick} precisely when it is closed under these
operations. Thus one is naturally led to the study  the lattice of thick
subcategories of $\F$.

The lattice of thick subcategories of $\F$ is determined by the Morava K-theories.
For each $n \ge 1$ there is a spectrum $K(n)$ called the $n$-th
Morava K-theory whose coefficient ring $K(n)_*$ is isomorphic to
$\mathbb{F}_p[v_n , v_n^{-1}]$ with $|v_n| = 2(p^n -1)$. We also set
$K(0)$ to be the rational Eilenberg-Mac Lane spectrum and $K(\infty)$
the mod-$p$ Eilenberg-Mac Lane spectrum. These theories have the
following pleasant properties.

\begin{enumerate}
\item For every spectrum $X$, $K(n)\wedge X$ has the homotopy type
of a wedge of suspensions of $K(n)$.
\item K\"unneth isomorphism: $K(n)_* (X \wedge Y) \cong K(n)_* X \otimes _{K(n)_*} K(n)_*
Y$. In particular $K(n)_* (X \wedge Y) = 0 $ if and only if either
$K(n)_* X = 0$ or $K(n)_* Y = 0$.
\item If $X \ne 0$ and finite, then for all $n > > 0$, $K(n)_* X  \ne 0$.
\item For each $n$, $K(n+1)_* X = 0$ implies $K(n)_* X = 0$
\item (Nilpotence theorem) Morava K-theories detect ring spectra: If
$R$ is a non-trivial ring spectrum, then there exists an $n$ ($0 \le
n \le \infty$) such that $K(n)_* R \ne 0$
\end{enumerate}

The first three properties can be easily derived from the fact that
every graded module over $K(n)_*$ is a direct sum of suspensions of $K(n)$.
The third property is proved in \cite{rav}, and the last property can be derived from the  $MU$-version of the
nilpotence theorem stated above; see \cite{hs} for a proof of this implication. In view of the above properties,
one prefers to work with $K(n)$ as opposed to $MU$ because it is
easier to do computations with $K(n)$. Set $\C_0 = \F$, and for $n
\ge 1$, let  $\C_n := \{ X \in \F : K(n-1)_* X = 0\}$, and finally
let  $\C_{\infty}$  denote  the subcategory of contractible spectra.
We can now state the celebrated thick subcategory theorem.

\begin{thm}\cite{hs} \emph{(Thick subcategory theorem)} A subcategory $\C$ of $\F$ is thick if and only if $\C = \C_n$
for some $n$. Further these subcategories form a nested decreasing
filtration of $\F$:
\[ \C_{\infty} \subsetneq \cdots \subsetneq \C_{n+1} \subsetneq \C_n \subsetneq \C_{n-1} \subsetneq \cdots
\subsetneq \C_1 \subsetneq \C_0 \]
\end{thm}

We say that a spectrum $X$ is of type-$n$ if it belongs to $\C_n -
\C_{n+1}$, and we write $\text{type}(X) = n$.  For example the sphere spectrum is of type $0$ and the
mod-p Moore spectrum is of type-$1$. The existence of type-$n$ spectra was first proved by Mitchell \cite{mit}.

It is not hard to prove the above theorem using the nilpotence
theorem. The proof we sketch below is following Rickard.  The only
other tool that we need is finite localisation.

\begin{thm} \cite{miller} \emph{(Finite Localisation)}  Let $\C$ be a thick subcategory of $\F$, and let $\D$ denote
the localising subcategory generated by $\C$. Then there is a
localisation functor $\mathcal{L}^f_{\C}: \mathcal{S} \rar
\mathcal{S}$ called "finite localisation away from $\C$" which has
the following properties.
\begin{enumerate}
\item For $X$ finite, $\mathcal{L}^f_{\C} X = 0$  if and only if $X$ belongs
to $\F$.
\item For $X$ arbitrary, $\mathcal{L}^f_{\C} X = 0$  if and only if $X$ belongs
to $\D$.
\item $\mathcal{L}^f_{\C}$ is a smashing localisation functor, i.e., $\mathcal{L}^f_{\C} X \cong \mathcal{L}^f_{\C} S^0 \wedge X$, where $S^0$
is the $p$-local sphere spectrum.
\end{enumerate}
\end{thm}

The idea involved in the construction of such a localisation functor is well-known to homotopy theorists. For instance, it shows up in the proof  of the
Brown representability theorem. A very good treatment of this construction is also given by Rickard \cite{Ri} where he constructs idempotent modules in the
stable module category using finite localisation.

We should mention at this point that a version of Ravenel's
telescope conjecture states that every smashing localisation functor
(a localisation functor that satisfies property (3) above) on
$\mathcal{S}$ is isomorphic to $\mathcal{L}^f_{\C_n}$ for some integer $n$.
This is the only conjecture of Ravenel that is still open;
some experts believe that it is false, see \cite{rav-tc}.

Now we give a proof (due to Rickard) of the thick subcategory
theorem. Let $\C$ be a non-zero thick subcategory of $\F$. Then
define
\[ n := \text{max}\; \{l: \C \subseteq \C_l \}.\]
From property (3) of the Morava K-theories we infer that $n$ is a
well-defined non-negative integer. We claim that $\C = \C_n$. Note
that we only have to show that $\C_n \subseteq \C$. In showing this
inclusion we use the finite localisation functors
$\mathcal{L}^f_{\C}$. So let $X$ in $\C_n$. Now to show that $X$ is
in $\C$, it is enough to show that $\mathcal{L}^f_{\C}\, X = 0$. But
since these functors are smashing we have to show $X \wedge
\mathcal{L}^f_{\C} \, S^0 = 0$. Since every finite spectrum $(X)$ is
Bousfield equivalent to a ring spectrum $(X \wedge DX)$, we can
assume without loss of generality that $X$ is a ring spectrum. Then
by property (5) it suffices to show that for all $0 \le l \le
\infty$, $K(l)_* \,(X \wedge \mathcal{L}^f_{\C} S^0) = 0$.  Further
by property (2) we have to show that for each $l$, either $K(l)_*
\,X = 0$ or $K(l)_*(\mathcal{L}^f_{\C} S^0) = 0$. Since $X$ is in
$\C_n$, the former holds for all $0 \le l < n$ by property (4). So we have
to show that the latter holds for $n \le l \le \infty$. Now by the
definition of $n$, we have for each $n \le l \le \infty$, a spectrum $X_l$ in
$\C$ such that $K(l)_*(X_l) \ne 0$. Since $X_l$ is in $\C$, we have
$\mathcal{L}^f_{\C} X_l = X_l \wedge \mathcal{L}^f_{\C}S^0 = 0$. So
clearly $K(l)_* \, \mathcal{L}^f_{\C} X_l = 0$, but since $K(l)_*\,
(X_l) \ne 0$, we must have for all $n \le l \le \infty$, that $K(l)_* \,
\mathcal{L}^f_{\C} S^0 = 0$ as desired.

Note that this proof highlights the key properties of Morava
K-theories which are used in proving the thick subcategory theorem, and
therefore it can be adapted easily to the other algebraic settings
such as derived categories of rings and stable module categories of
group algebras. The role played by the Morava K-theories in the
former are the residue fields and in the latter are the kappa
modules; see \cite{mps} for a thick subcategory theorem in an
axiomatic stable homotopy category.

We now illustrate how one can use the thick subcategory theorem.
Suppose P is some generic property of spectra and we want to
identify the subcategory of finite spectra which satisfy P. If we can find a type-$k$ spectrum which satisfies P and
a type-$(k-1)$ spectrum which does not satisfy P, that forces the
subcategory in question to be $\C_k$. For example, consider the
generalised homology theory $BP \la n \ra$ whose coefficient ring is
given by $\mathbb{Z}_p[v_1, v_2, \cdots v_n]$ with $|v_i| = 2(p^i -
1)$.  Using the above strategy one can easily show that the full
subcategory of finite spectra which have bounded $BP \la n \ra$
homology  (spectra $X$ such that $BP \la n \ra_i X = 0$ for $i > >
0$) is precisely $\C_{n+1}$.

\section{Bousfield classes of finite spectra}

There are several interesting applications of the thick subcategory
theorem. We focus on its applications to the Bousfield lattice -- an
important lattice which encapsulates  the gross structure of
stable homotopy theory. This was introduced by Bousfield in
\cite{bob, boa}. Given a spectrum $E$, define its \emph{Bousfield Class} $\la
E \ra$ to be the collection of all spectra which are invisible to
the $E$-homology theory, i.e., spectra $X$ such that $E_*(X) = 0$ or
equivalently $E \wedge X = 0$. Then we say that spectra $E$ and $F$
are Bousfield equivalent if $\la E \ra = \la F \ra.$ It is a result
by Ohkawa that there is only a set of Bousfield classes. With the
partial  order given by reverse inclusion, the set of Bousfield
classes form a lattice which is called the Bousfield lattice and will
be denoted by $\mathbf{B}$.  One can perform various operations on
$\mathbf{B}$. The two important ones being the wedge $(\vee)$ : $\la
X \ra \vee \la Y \ra = \la X \vee Y \ra $ and the smash $(\wedge)$:
 $\la X \ra \wedge \la Y \ra = \la X \wedge Y \ra $. In this
lattice, the Bousfield class of the sphere spectrum is the largest element
and that of the trivial spectrum is the smallest. This lattice plays
an important role in the study of modern stable homotopy theory.
While much of the current knowledge about the $\mathbf{B}$ is only
conjectural, the thick subcategory theorem completes determines the
Bousfield classes of finite spectra. We describe them in the next
theorem which was conjectured by Ravenel.

\begin{thm}\cite{hs} \emph{(Class-invariance theorem)} Let $X$ and $Y$ be finite $p$-local spectra, then
$\la X \ra \le \la Y \ra$ if and only if $\text{type}(X) \ge \text{type}(Y)$.
\end{thm}

Although this theorem follows as an immediate corollary to the thick subcategory theorem, it is a very non-trivial statement
about finite spectra. It says that the Bousfield class of a finite spectrum is completely determined by its type.

We now discuss the Boolean algebra conjecture of Ravenel which
identifies the Boolean subalgebra generated by the Bousfield classes
of finite $p$-local spectra. But first we have to introduce some important
non-nilpotent maps of finite spectra called $v_n$-self maps. A self
map $f: \Sigma^? X \rar X$ is a $v_n$-self map ($n \ge 1$) if
$K(n)_*(f)$ is an isomorphism and $K(m)_*(f)$ is zero for $m \ne n$.
For example, the degree $p$ map on the sphere spectrum is a $v_0$-self map, and the
Adams map \cite{adamsmap} on the Moore spectrum: $\Sigma^? M(p) \rar M(p)$ which
induces isomorphism in complex $K$-theory is a $v_1$-self map. These
$v_n$-self maps are important because give rise to periodic families
in the stable homotopy groups of spheres. For example, one can
iterate Adams maps and get a periodic family in $\pi_*(S^0)$ called the
$\alpha$-family.
Showing the existence of such maps is highly non-trivial.
A deep result of Hopkins and Smith called the periodicity
theorem produces a wealth of such maps. More precisely:

\begin{thm}\cite{hs} \emph{(Periodicity Theorem)}
\begin{enumerate}
\item Every type-n spectrum admits an asymptotically unique $v_n$-self map $\phi_X: \Sigma^?
X \rar X$
\item If $h: X \rar Y$ is a map between type-n spectra,
then there exits integers $i$ and $j$ such that the follow diagram
commutes:
\[
\xymatrix{
\Sigma^? X \ar[r]^{\Sigma^? h} \ar[d]_{\phi_X^i} & \Sigma^? Y \ar[d]^{\phi_Y^j} \\
X \ar[r]^{h} & Y
}
\]
 \end{enumerate}
\end{thm}

Using this periodicity theorem it is not hard to show that the full
subcategory of finite $p$-local spectra admitting $v_n$ self map is
precisely $\C_n$. So this theorem gives another characterisation of the thick subcategories of $\F$.

For every positive integer $n$, let $F(n)$ denote some
spectrum of type-$n$. Note that the Bousfield class of $F(n)$ is well-defined by the class-invariance theorem.
Now the periodicity theorem says that $F(n)$ admits an essentially unique $v_n$-self map. So let $T(n)$ denote
the mapping telescope of this $v_n$-self map.  It follows that
the Bousfield classes of $T(n)$ is also  well-defined.

A Bousfield class $\la E \ra$ is said to be complemented if there
exists another class $\la F \ra$ such that $\la E \ra \wedge \la F
\ra = \la 0 \ra$ and $\la E \ra \vee \la F \ra = \la S^0 \ra$. The
collection of all complemented Bousfield classes forms a Boolean
algebra with respect to the smash and wedge operations and will be
denoted by $\mathbf{BA}$.
Bousfield \cite{boa} showed that every possibly infinite wedge of
finite spectra belongs to $\mathbf{BA}$. A  pleasant
consequence of the thick subcategory theorem is the classification
of the Boolean subalgebra generated by the finite $p$-local spectra
and their complements in the $p$-local sphere spectrum.

\begin{thm} \cite{rav} \emph{(Boolean Algebra Theorem)} Let $\mathbf{FBA}$
denote the Boolean subalgebra generated by the Bousfield classes of the finite $p$-local
spectra and their complements in $\la S^0 \ra$. Then $\mathbf{FBA}$ is the free (under
complements, finite unions and finite intersections) Boolean algebra
generated  by the Bousfield classes of the telescopes $\la T(n) \ra$
for $n\ge 0$.
\end{thm}

So by this theorem one can identify $\mathbf{FBA}$ with the Boolean algebra of finite  and cofinite subsets of non-negative integers:
the Bousfield class $\la T(n) \ra$ corresponds to the subset   $\{ n \}$, and $\la F(n) \ra$ corresponds to the subset $\{ n, n+1, n+2, \cdots \}$.
By the way, the Boolean algebra conjecture of Ravenel uses $K(n)$ instead of
$T(n)$ in the above theorem; according to telescope conjecture
(which is still open) these two spectra are Bousfield equivalent.

There are several other interesting sublattices of $\mathbf{B}$ which have been studied.
For example there is a distributive lattice $\mathbf{DL}$ which consists of the Bousfield
classes $\la X \ra$ such that $\la X \ra \wedge \la X \ra = \la X \ra$ which has some nice properties.
A good discussion on the structure of the Bousfield lattice can be found in \cite{hovpallattice}. These
authors use lattice theoretic methods to explore the structure of the Bousfield lattice. They also pose
a number of interesting conjectures and study their implications.

We end by mentioning briefly one other application of the thick
subcategory theorem.  Thomason has given a brilliant K-theory recipe
\cite{Th} which  refines the thick subcategory theorem and gives a classification of the triangulated subcategories of finite
spectra. This recipe amounts to computing the Grothendieck groups of the
thick subcategories of the finite $p$-local spectra. We refer the
reader to \cite{Cheb06a} where we use this recipe to study the
lattice of triangulated subcategories of finite spectra.

\bibliographystyle{alpha}

\end{document}